\documentstyle[12pt]{article}
\newcommand{\nz}{{\rm I\!N}}

\newcommand{\rz}{{\rm I\!R}}
\newcommand{\cz}{{\,\rm{\sf I}\!\!\!C}}

\begin{document}
\begin{center}
{\bf LIMIT THEOREMS FOR THE HIERARCHY OF FREENESS}
\\[40pt]
{\sc Uwe Franz}$^{\dagger}$ {\sc and Romuald Lenczewski}$^{\ddagger}$
\footnote{ This work is supported by KBN (grant No. 2PO3A 05415).}
\\[40pt]
$^{\dagger}$ Institut de Recherche Math{\'e}matique Avanc{\'e}e,\\
Universit$\acute{{\rm e}}$ Louis Pasteur,\\
Rue Ren{\'e} Descartes,\\
F-67084 Strasbourg Cedex, France\\
e-mail franz@math.u-strasbg.fr\\[20pt]
$^{\ddagger}$ Hugo Steinhaus Center for Stochastic Methods,\\
Institute of Mathematics,\\
Technical University of Wroclaw,\\
Wybrze$\dot{{\rm z}}$e Wyspianskiego 27\\
50-370 Wroclaw, Poland\\
e-mail lenczew@im.pwr.wroc.pl\\[20pt]
\end{center}
\begin{abstract}
The central limit theorem, the invariance principle and
the Poisson limit theorem for the hierarchy of freeness 
are studied. We show that for given
$m\in \nz$ the limit laws can be expressed in terms of
non-crossing partitions of depth smaller or equal to $m$.
For ${\cal A}=\cz [x]$, we solve the associated 
moment problems and find explicitly 
the discrete limit measures.
\end{abstract}
\vspace{20pt}
\begin{center}
{\sc 1. Introduction}
\end{center}
The notion of the {\it hierarchy of freeness} was introduced
in [Len97] in the context of a unification of the main types of
non-commutative independence (tensor, free, and Boolean, see
the axiomatic approach in [Sch94,Sch95]).
The main idea of the construction presented in [Len97] was to
approximate the free product of states [Voi85] through a sequence of
products called $m$-{\it free products}, $m\in \nz$, using only tensor
independence. 
In this way one obtains a hierarchy of products as well as a
hierarchy of non-commutative probability spaces, of which the
latter was called in [Len97] the hierarchy of freeness.
 
In the hierarchy of $m$-free products the two extremes are given by
the Boolean product which corresponds to the first order
approximation $m=1$ and the free product, obtained for
$m=\infty$. Thus the hierarchy fills the ``gap'' between the
Boolean product and the free product. Its another important feature
is that it equips the combinatorics of non-crossing partitions
with a hierarchic structure induced by their depths. 
Recall that the combinatorics of the Boolean 
product is based on the so-called {\it interval partitions} and
that of the free product 
-- on all {\it non-crossing partitions}. 
By studying convolution-type limit theorems in this paper, 
we establish a connection between the combinatorics of the
$m$-free product (or, rather of the $m$-free convolution) 
and non-crossing partitions of depth
$d(P)\leq m$. Thus the hierarchy also fills the ``gap'' between
the combinatorics of interval partitions and that of 
all non-crossing partitions. 
Let us add that the hierarchy of freeness lends itself
easily to certain generalizations, and in fact was introduced in
[Len97] in the context of the conditionally free product [BLS96]
of states. Other generalizations were indicated in [FLS98].

In this work we study the convolution-type central limit
theorems, the invariance principles and Poisson's limit theorems
for $m$-free products, calling those theorems $m$-free limit
theorems. Let us only note that we do not use the $m$-free
convolutions in our notations. Nevertheless, all theorems can be phrased 
using $m$-free convolutions introduced in [Len97]. 
It is well known that in the central limit theorem for free
independence [Voi85] only non-crossing pair partitions give rise
to the limit Wigner semi-circle law [Spe90].
In our case we show that in the $m$-free central limit theorem 
only non-crossing pair partitions of depth less
than or equal to $m$ appear in the combinatorial form of the 
limit law for each $m\in \nz$.
For the special case of the algebra of polynomials in
one variable $\cz[x]$, we introduce a hierarchy of Cauchy
transforms of the limit laws which enables us to recover
the corresponding hierarchy of discrete measures on the real
line which approximate the Wigner measure. A similar approach is
used for $m$-free Poisson's limit theorems.

Section 2 is of preliminary character and contains all needed
facts on the hierarchy of freeness. In Section 3 we prove the
central limit theorem for the hierarchy of freeness (Theorem
3.5). Note that our approach is based on the tensor product
construction developed in [Len97] and as such gives a new
(and probably the most explicit) proof of the free central limit theorem.
In Section 4, the corresponding invariance principle is stated
(Theorem 4.1) and a hierarchy of $m$-free Brownian motions is introduced.
In Section 5, we restrict ourselves to $\cz [x]$ and
study the hierarchy of measures corresponding to the central
limit laws. We show that they are discrete measures that
approximate weakly the Wigner measure.
Poisson's limit theorem for the hierarchy of freeness is proved
in Section 6 and the associated moment problems are solved.
\\[5pt]
\begin{center}
{\sc 2. The Hierarchy of Freeness}
\end{center}
This section is of preliminary character and contains all needed
facts on the hierarchy of freeness. For more details, see [Len97]
and [FLS98]. 

Let $({\cal A}_{l})_{l\in I}$ be a family of unital
*-algebras and let 
$(\phi_{l})_{l\in I}$ be the corresponding family of states. 
We assume that ${\cal A}_{l}={\cal A}_{l}^{0}\oplus {\bf
1}_{l}$, where ${\cal A}_{l}^{0} $ is a *-subalgebra of ${\cal
A}_{l}$, and in the free product $*_{l\in I}{\cal A}_{l}$ we identify units. 
Extend each ${\cal A}_{l}$ to 
$\widetilde{\cal A}_{l}={\cal A}_{l}*{\cz}(t_{l})$, where
${\cz}(t_{l})$ is the unital *-algebra generated by the projection
$t_{l}$.
Make $\widetilde{\cal A}_{l}$ into a *-algebra in the canonical fashion.
Finally, denote by 
$(\widetilde{\phi}_{l})_{l\in I}$ the Boolean extensions of
$(\widetilde{\phi}_{l})_{l\in I}$, i.e. states on 
$(\widetilde{\cal A}_{l})_{l\in I}$ given by 
$\widetilde{\phi}_{l}({\bf 1}_{l})=1$ and 
$$
\widetilde{\phi_{l}}(t_{l}^{r}a^{(1)}t_{l}\ldots 
t_{l}a^{(p)}t_{l}^{s})
=\phi_{l}(a^{(1)})\ldots \phi_{l}(a^{(p)})
$$
for $a^{(1)}, \ldots , a^{(p)} \in {\cal A}_{l}^{0}$, 
$r,s\in \{0,1\}$. For details, see [FLS98].

Consider the quantum probability space $({\cal B}, \widetilde{\Phi})$, where
$$
{\cal B}= \bigotimes_{l\in I}\widetilde{\cal A}_{l}^{\otimes \infty},\;\;\;
\widetilde{\Phi}=\bigotimes_{l\in I}\widetilde{\phi}_{l}^{\otimes \infty},
$$
and the tensor products are understood as in 
[FLS98], with
canonical involutions on $\bigotimes_{l\in I}{\cal A}_{l}$ and ${\cal
B}$. This is the quantum probability space in which one can embed
the hierarchy of freeness defined in [Len97] (see again [FLS98]). 
Since we have two tensor
products here (over $I$ and then over ${\nz}$ for each $l\in
I$), we will label tensor sites by $(l,k)$, $l\in I$, $k\in {\nz}$ and we will refer to $l$ and $k$ as the outer and inner
site, respectively.

In the definition of these embeddings the following notations will be used.
For $l\in I, n\in {\nz}$, let
$$
i_{n}^{(l)}: \;\;\widetilde{\cal A}_{l}\rightarrow \widetilde{\cal
A}_{l}^{\otimes \infty}
$$
be the linear mapping given by
$$
i_{n}^{(l)}(a)={\bf 1}_{l}^{\otimes (n-1)}\otimes a \otimes
{\bf 1}_{l}^{\otimes \infty},
$$
for $a\in \widetilde{\cal A}_{l}$. For notational convenience we put 
$i_{0}^{(l)}(a)=0$. Further, we denote by
$$
t_{[k}^{(l)}={\bf 1}_{l}^{\otimes (k-1)}\otimes t_{l}^{\otimes \infty}
$$
a projection in $\widetilde{\cal A}_{l}^{\otimes \infty}$ which
is built from projections $t_{l}$ at all sites $\geq k$, $k\geq
1$, and we put for convenience  $t_{[0}^{(l)}=0$.

We define the linear mappings
$$
\gamma_{k}^{(l)}\;:\;\; {\cal A}_{l}^{0}\rightarrow {\cal B},\;\;\;\;
\gamma_{k}^{(l)}(a)=i_{k}^{(l)}(a)\otimes \bigotimes_{r\neq l}
t_{[k}^{(r)},
$$
$$
\widehat{\gamma}_{k}^{(l)}\;:\;\; {\cal A}_{l}^{0}\rightarrow {\cal B},\;\;\;\;
\widehat{\gamma}_{k}^{(l)}(a)=i_{k}^{(l)}(a)\otimes \bigotimes_{r\neq l}
t_{[k-1}^{(r)},
$$
where $k\in {\nz}$, $l\in I$. 
Note that since $i_{0}^{(l)}(t)=0$, we have 
$\widehat{\gamma}_{1}^{(l)}(a)=0$. 
In other words, $\gamma_{k}^{(l)}(a)$ puts $a\in {\cal
A}_{l}^{0}$ at site $(l,k)$
and projections $t_{r}$ at sites $(r,s)$ for all $r\neq l$ and 
$s\geq k$. In turn, $\widehat{\gamma}_{k}^{(l)}(a)$ puts $a$ at site
$(l,k)$ and projections $t_{r}$ at sites $(r,s)$ for all $r\neq
l$ and $s\geq k-1$.

It was shown in [FLS98] that the mappings
$$
j_{l}^{(m)}\;:\;\; {\cal A}_{l}^{0}\rightarrow {\cal B},
$$
$$
j_{l}^{(m)}=\sum_{k=1}^{m}j_{l,k}\equiv 
\sum_{k=1}^{m}(\gamma_{k}^{(l)}-\widehat{\gamma}_{k}^{(l)}),
$$
where $l\in I$, $m\in {\nz}$, are *-homomorphisms. Using them,
we can define for each $m\in {\nz}$ the *-homomorphism
$$
j^{(m)}: *_{l\in I}{\cal A}_{l}\rightarrow {\cal B}
$$
as the linear extension of
$j^{(m)}({\bf 1})=\bigotimes_{l\in I}{\bf 1}_{l}^{\otimes \infty}$ and
$$
j^{(m)}(a_{1}\ldots a_{n})=j_{l_{1}}^{(m)}(a_{1})\ldots 
j_{l_{n}}^{(m)}(a_{n}),
$$
where $a_{i} \in {\cal A}_{l_{i}}^{0}$, $l_{i}\in I$, $i=1, \ldots , n$.
\\[0pt]
\indent{\par}
{\sc Definition 2.1} 
The sequence of quantum probability spaces
$({\cal A}^{(m)}, \Phi^{(m)})_{m\in {\nz}}$, where
${\cal A}^{(m)}=j^{(m)}(*_{l\in I}{\cal A}_{l})$ and 
$\Phi^{(m)}$ is the restriction of $\widetilde{\Phi}$ to
${\cal A}^{(m)}$, is called the {\it hierarchy of freeness}.
The state $\widetilde{\Phi}^{(m)}$ is called the 
$m$-free product state and
$j^{(m)}(a)$, $a\in {\cal A}_{l}^{0}$ are called the 
$m$-{\it free random variables}.\\[0pt]
\indent{\par}
{\sc Remark.}
Note that $\widetilde{\Phi}\circ j^{(m)}$ defines a state on
$*_{l\in I}{\cal A}_{l}$.\\[0pt]
\indent{\par}
The GNS construction for the hierarchy of freeness 
[FLS98] will also be useful here.
Thus, let (${\cal H}_{l}, \pi_{l}, \Omega_{l}$) be the GNS triple
associated with the pair (${\cal A}_{l}, \phi_{l}$), i.e. ${\cal
H}_{l}$ is a pre-Hilbert space, $\pi_{l}$ is a *-representation
of ${\cal A}_{l}$ and $\Omega_{l}$ is a cyclic vector, such that
$
\phi_{l}(x)=\langle \Omega_{l}, \pi_{l}(x)\Omega_{l} \rangle
$
for any $x\in {\cal A}_{l}$.
We start from the infinite tensor product pre-Hilbert space 
$$
{\cal H}^{\otimes}=\bigotimes_{l\in I}{\cal H}_{l}^{\otimes \infty}
$$
with respect to the vector 
$\Omega=\bigotimes_{l\in I}\Omega_{l}^{\otimes \infty}$
and denote by
$$
\Gamma_{k}^{(l)}\;:\;\; 
{\cal A}_{l}^{0}\rightarrow {\cal L}({\cal H}^{\otimes}),\;\;\;
\widehat{\Gamma}_{k}^{(l)}\;:\;\;
{\cal A}_{l}^{0}\rightarrow {\cal L}({\cal H}^{\otimes})
$$
the *-homomorphisms corresponding to $\gamma_{k}^{(l)}$, 
$\widehat{\gamma}_{k}^{(l)}$, i.e.
$$
\Gamma_{k}^{(l)}(a)=
i_{k}^{(l)}(\pi_{l}(a))\otimes
\bigotimes_{j\neq l}
P_{[k}^{(j)},
$$
$$
\widehat{\Gamma}_{k}^{(l)}(a)=
i_{k}^{(l)}(\pi_{l}(a))\otimes
\bigotimes_{j\neq l}
P_{[k-1}^{(j)}
$$
for $a\in {\cal A}_{l}^{0}$, where
$P_{[k}^{(j)}={\rm Id}^{\otimes (k-1)}\otimes (P^{(j)})^{\otimes
\infty}$, $P^{(j)}$ is the projection onto the vacuum $\Omega_j$ in
${\cal H}_j$, and $P_{[0}=0$.
Then the GNS representation $\pi^{\otimes m}$ of 
$(*_{l\in I}{\cal A}_{l}, \Phi\circ j^{(m)})$ 
is given by $\pi^{\otimes m}({\bf 1})=
\bigotimes_{l\in I}{\rm Id}_{l}^{\otimes \infty}$ and
$\pi^{\otimes m}=*_{l\in I}\pi_{l}^{\otimes m}$ on 
$*_{l\in I}{\cal A}_{l}^{0}$, where
$$
\pi_{l}^{\otimes m}(a)=
\sum_{k=1}^{m}
(\Gamma_{k}^{(l)}(a)-\widehat{\Gamma}_{k}^{(l)}(a))
$$
for $a\in {\cal A}_{l}^{0}$.
For each $m\in {\nz}$ the cyclic vector is $\Omega$
and the carrier space of $\pi^{\otimes m}$ is
${\cal H}^{\otimes m}=\pi^{\otimes m}(*_{l \in I}{\cal A}_{l})
\Omega.$

We need to take a closer look at the correlations
$$
\widetilde{\Phi}
\left(
j_{l_{1}}^{(m)}(a_{1})\ldots j_{l_{n}}^{(m)}(a_{n})
\right)
=\sum_{1\leq m_{1}, \ldots , m_{n} \leq m}
\widetilde{\Phi}
\left(
j_{l_{1},m_{1}}(a_{1})
\ldots 
j_{l_{n}, m_{n}}(a_{n})
\right)
$$
for any tuple $(l_{1}, \ldots , l_{n})$, $a_{i}\in {\cal
A}_{l_{i}}^{0}$, $i=1, \ldots , n$. Equivalently, we can write
$$
\widetilde{\Phi}
\left(
j_{l_{1}}^{(m)}(a_{1})\ldots j_{l_{n}}^{(m)}(a_{n})
\right)=
\widetilde{\Phi}
\circ j^{(m)}(a_{1} \ldots , a_{n}).
$$
Before we derive some 
results which are specific to the central limit theorem and use
the assumption on the zero mean, we
prove a ``pyramid formula''(slightly more general than the one
in [Len97] which always allows us to reduce the summation in the
above sum to a ``pyramid''. We also give a new proof, using
the GNS construction.\\[0pt]
\indent{\par}
{\sc Proposition 2.2}
{\it The following formula holds:}
$$
\widetilde{\Phi}
\left(
j_{l_{1}}^{(m)}(a_{1})\ldots j_{l_{n}}^{(m)}(a_{n})
\right)
=\sum_{(m_{1}, \ldots , m_{n}) \in \Upsilon _{n}^{m}}
\widetilde{\Phi}
\left(
j_{l_{1},m_{1}}(a_{1})
\ldots 
j_{l_{n}, m_{n}}(a_{n})
\right),
$$
{\it where
$\Upsilon _{n}^{m}=\{(p_{1}, \ldots , p_{n})| 1\leq p_{k},
p_{n-k} \leq k\wedge m, 1\leq k \leq n/2\}$ 
and $k\wedge m={\rm min}\{k,m\}$}.\\[5pt]
{\it Proof.}
Using the GNS construction, we obtain
$$
\widetilde{\Phi}\left(
j_{l_{1}}^{(m)}(a_{1})\ldots j_{l_{n}}^{(m)}(a_{n})
\right)
=\langle \Omega , 
\pi^{\otimes m}(a_{1})\ldots , \pi^{\otimes m}(a_{n})
\Omega
\rangle
$$
and thus, in order to prove the proposition, it is enough
to show that if $(m_{1}, \ldots , m_{n})\notin \Upsilon_{n}^{m}$,
then 
$$
\langle\Omega,
(\Gamma_{m_{1}}^{(l_{1})}(a_{1})-\widehat{\Gamma}_{m_{1}}^{(l_{1})}(a_{1}))
\ldots
(\Gamma_{m_{n}}^{(l_{n})}(a_{n})-\widehat{\Gamma}_{m_{n}}^{(l_{n})}(a_{n}))
\Omega\rangle =0.
$$
Introduce the filtration 
$$
{\cal H}^{\otimes}_{0]} \subset {\cal H}^{\otimes}_{1]}
\subset \ldots \subset {\cal H}^{\otimes}_{k]}\subset \ldots 
$$
of subspaces of ${\cal H}^{\otimes}$ given by ${\cal
H}_{0]}={\cz}\Omega$ and
$$
{\cal H}^{\otimes}_{k]}={\rm Lin}\;\{\bigotimes_{l\in
I}(x_{l,1}\otimes \ldots \otimes x_{l,k}\otimes
\Omega_{l}^{\otimes \infty})\}.
$$
Note that if $k>1$, then 
$\Gamma_{k}^{(l)}(a)$ agrees with $\widehat{\Gamma}_{k}^{(l)}(a)$ on
${\cal H}_{k-2]}$. Moreover, 
$$
\left(
\Gamma_{k}^{(l)}(a)-\widehat{\Gamma}_{k}^{(l)}(a)
\right)\;
{\cal H}^{\otimes}_{k-1]}\subset {\cal H}^{\otimes}_{k]},
$$
for any $k\geq 1$. These two facts imply that we have
$$
(\Gamma_{m_{1}}^{(l_{1})}(a_{1})-\widehat{\Gamma}_{m_{1}}^{(l_{1})}(a_{1}))
\ldots
(\Gamma_{m_{n}}^{(l_{n})}(a_{n})-\widehat{\Gamma}_{m_{n}}^{(l_{n})}(a_{n}))
\Omega =0
$$
if $(m_{1}, \ldots , m_{n})\notin \Theta _{n}^{m}$, where
$$
\Theta _{n}^{m}=\{(p_{1}, \ldots , p_{n})| 1\leq p_{i} \leq
(n-i+1)\wedge m\}.
$$ 
We can repeat this argument for the adjoints and obtain a mirror
reflection of this condition 
($(m_{n}, \ldots , m_{1})\notin \Theta _{n}^{m}$), 
which finally leads to 
$$
\langle \Omega ,
(\Gamma_{m_{1}}^{(l_{1})}(a_{1})-\widehat{\Gamma}_{m_{1}}^{(l_{1})}(a_{1}))
\ldots
(\Gamma_{m_{n}}^{(l_{n})}(a_{n})-\widehat{\Gamma}_{m_{n}}^{(l_{n})}(a_{n}))
\Omega\rangle =0
$$
if $(m_{1}, \ldots , m_{n})\notin \Upsilon _{n}^{m}$. 
\hfill $\Box$\\[0pt]
\indent{\par}
{\sc Proposition 2.3.}
{\it If ${\cal A}_{l}={\cal A}$, $\phi_{l}=\phi$, $l\in I$, then
the correlations of $m$-free random variables 
are invariant under permutations $\pi$ of ${\nz}$, i.e.}
$$
\widetilde{\Phi}
\left(
j^{(m)}_{\pi(l_{1})}(a_{1})\ldots j^{(m)}_{\pi(l_{n})}(a_{n})
\right)
=\widetilde{\Phi}
\left(
j^{(m)}_{l_{1}}(a_{1})\ldots j^{(m)}_{l_{n}}(a_{n})
\right).
$$
{\it Moreover, if 
$\{l_{1}, \ldots , l_{r}\}\cap \{l_{r+1}, \ldots ,l_{n}\}=\emptyset$, then} 
$$
\widetilde{\Phi}(j^{(m)}_{l_{1}}(a_{1})\ldots j^{(m)}_{l_{n}}(a_{n}))=
\widetilde{\Phi}(j^{(m)}_{l_{1}}(a_{1})\ldots j^{(m)}_{l_{r}}(a_{r}))
\widetilde{\Phi}(j^{(m)}_{l_{r+1}}(a_{r+1})\ldots j^{(m)}_{l_{n}}(a_{n})).
$$
{\it Proof.}
From the properties of the tensor product and the fact that
$\phi_{l}=\phi$ for all $l\in I$, we obtain
$$
\widetilde{\Phi}
\left(
\gamma_{k_{1}}^{\natural (\pi(l_{1}))}(a_{1})
\ldots
\gamma_{k_{n}}^{\natural (\pi(l_{n}))}(a_{n})
\right)=
\widetilde{\Phi}
\left(
\gamma_{k_{1}}^{\natural (l_{1})}(a_{1})
\ldots
\gamma_{k_{n}}^{\natural (l_{n})}(a_{n})
\right)
$$
for any $1\leq k_{1}, \ldots , k_{n} \leq m$, where
$\gamma_{k}^{\natural (l)}(a)=\gamma_{k}^{(l)}(a), 
\widehat{\gamma}_{k}^{(l)}(a)$. From this follows the first
part of the proposition. The second part is obvious. 
\hfill $\Box$\\[5pt]
\begin{center}
{\sc 3. A Central Limit Theorem}
\end{center}
In this section we prove the central limit theorem for the sums
of $m$-free independent random variables. We show that in the limit
only the non-crossing pair partitions $P$ of depth $d(P)\leq m$ give a
nonvanishing contribution.\\[0pt]
\indent{\par}
{\sc Definition 3.1.}
A pair partition 
$P=\{P_{1}, \ldots P_{k}\}$, 
where $P_{j}=\{\alpha (j),\beta (j)\}$, $j=1, \ldots , k$, 
of the set $\{1, \ldots , 2k\}$ 
is {\it crossing} if there exist $1\leq p,q\leq k$
such that $\alpha(p)<\alpha(q)<\beta(p)<\beta(q)$. If $P$ is not a crossing
partition, then it is called {\it non-crossing}.
If $P$ is non-crossing, then by $d(P)$ we denote its 
{\it depth}, i.e. the maximal of all integers $d$, 
for which there exist $1\leq s_{1}, \ldots , s_{d}\leq k$ such that
$\alpha(s_{1})<\ldots <\alpha(s_{d})$ and 
$\beta(s_{1})> \ldots > \beta(s_{d})$. We will denote the set of
all non-crossing pair partitions $P$ of depth $d(P)\le m$ of the set
$\{1,\ldots,n \}$ by $NC_{n}^{{\rm pair}}(m)$.\\[0pt]
\indent{\par}
{\sc Remark.} If we link each  $\alpha(l)$ with $\beta(l)$ in 
a pair-partition $P$ by ``bridges'', then a pair
partition is non-crossing if and only if it is possible to draw these
bridges without intersections. The depth $d(P)$ of $P$ is then
the maximal number of bridges that pass over the same ``gap''.

Note that with each tuple 
$(l_{1}, \ldots , l_{n})$, $l_{1}, \ldots , l_{n}\in I$,
we can associate a partition $P$ of $\{1, \ldots , n\}$. 
This can be done as follows.
Let $K=\{k_{1}, \ldots , k_{r}\}=\{l_{1}, \ldots , l_{n}\}$
with $k_{1}<k_{2} <\ldots < k_{r}$ and put 
$$
P_{i}=\{p|\;k_{p}=i\}.
$$ 
Then we will say that the partition $P$ is associated with the tuple
$(l_{1}, \ldots , l_{n})$.\\[0pt]
\indent{\par}
{\sc Lemma 3.2.}
{\it Assume that the partition $P$ 
associated with the tuple $(l_{1}, \ldots , l_{n})$, where $n=2k$,
is a non-crossing pair-partition of depth $d(P)>m$. 
If $\phi(a_{i})=0$ for $i=1, \ldots , n$, then}
$$
\widetilde{\Phi}
\left(
j^{(m)}_{l_{1}}(a_{1})\ldots j^{(m)}_{l_{n}}(a_{n})
\right)
=0
$$
{\it Proof.}
First of all note that each site can be occupied by at most two
elements since $P$ is a pair-partition. 
Assume that $d(P)>m$. Each $j_{l_{r}}^{(m)}(a_{r})$, $1\leq r
\leq n$ is a sum of $m$ terms in which $a_{r}$ appears at $m$ different 
sites, namely $(l_{r}, u)$, $1\leq u \leq m$. 
Since $P$ is a pair-partition and thus a given $a_{r}$ has only one
``partner'', say $a_{s}$ at site $(l_{s}, w)$ with
$l_{s}=l_{r}=l$, the only way to avoid ``singletons''
(first-order moments) is for each pair to occupy the same inner
site, i.e. $u=w$.
Now, we have at least $d(P)$ pairs to occupy at most 
$m$ different inner sites.
Since $d(P)>m$, at least one inner site, say $u$, must be
occupied by two pairs,
say $(a_{r}, a_{s})$ and $(a_{p}, a_{q})$, $l_{r}=l_{s}=l$,
$l_{p}=l_{q}=l'$.
Now, since $P$ is non-crossing, we 
must have $r<p<q<s$ or $p<r<s<q$. 
In the first case, at site $(l,u)$ we obtain
$$
\ldots a_{r}t\ldots ta_{s}\ldots 
$$
since $j_{l',u}(a_{p})$ and $j_{l',u}(a_{q})$
put a projection $t$ at all sites $(b,c)$, $b\neq l'$ and 
$c\geq u$. Thus $a_{r}$ and $a_{s}$ are separated by $t$ which
produces first moments, 
therefore gives zero by our zero mean assumption.
The second case is analogous.
\hfill $\Box$\\[0pt]
\indent{\par}
{\sc Lemma 3.3.}
{\it Assume that the partition $P$ 
associated with the tuple $(l_{1}, \ldots , l_{n})$, where
$n=2k$, is a non-crossing pair-partition of depth
$d(P)\leq m$. If $\phi(a_{i})=0$ for $i=1, \ldots , n$, then}
$$
\widetilde{\Phi}
\left(
j^{(m)}_{l_{1}}(a_{1})\ldots j^{(m)}_{l_{n}}(a_{n})
\right)
=\prod_{i=1}^{k}\phi(a_{P_{i}})
$$
{\it where $a_{J}=\prod_{l\in J} a_{l}$ for any $J\subset \{1,
\ldots , n\}$, with the product taken in the natural order.}\\[5pt]
{\it Proof.}
The proof will proceed by induction. 
Clearly, the case $m=1$ boils down to considering interval
pair-partitions (only they can be of depth $d(P)\leq 1$), i.e. 
take $P=\{\{i_{1},i_{2}\}, \ldots , \{i_{2k-1},i_{2k}\}\}$. Then
$$
\widetilde{\Phi}
\left(
j_{l_{1}}^{(1)}(a_{1})\ldots j_{l_{2k}}^{(1)}(a_{2k})
\right)
=\phi(a_{1}a_{2})\ldots \phi(a_{2k-1}a_{2k}).
$$
Assume now that 
$$
\widetilde{\Phi}
\left(
j^{(m-1)}_{l_{1}}(a_{1})\ldots j^{(m-1)}_{l_{n}}(a_{2k})
\right)
=
\prod_{i=1}^{k}\phi(a_{P_{i}})
$$
for $d(P)\leq m-1$ and any $k$. 
We will show that the same property holds
for $j^{(m)}$ and non-crossing partitions of depth
$d(P)\leq m$. 

The proof of that fact will be carried out by induction with
respect to $k$. If $k=1$, then we clearly have
$$
\widetilde{\Phi}
\left(
j_{l_{1}}^{(m)}(a_{1})j_{l_{2}}^{(m)}(a_{2})
\right)
=\phi(a_{1}a_{2}).
$$
Assume that 
$$
\widetilde{\Phi}
\left(
j^{(m)}_{l_{1}}(a_{1})\ldots j^{(m)}_{l_{2k-2}}(a_{2k-2})
\right)
=
\prod_{i=1}^{k-1}\phi(a_{S_{i}})
$$
for any tuple $(l_{1}, \ldots , l_{2k-2})$,
where $S$ is the partition associated with it
and $d(S)\leq m$. Now, when considering 
$\widetilde{\Phi}
\left(
j_{l_{1}}^{(m)}(a_{1})\ldots j_{l_{2k}}^{(m)}(a_{2k})
\right)$, it is enough to consider the case 
when $l_{1}=l_{2k}$ since
otherwise $P$ would separate into subpartitions and the
correlation would factorize by Proposition 2.3, thus we could
apply the inductive assumption with respect to $k$. 
By Proposition 2.2,
$$
\widetilde{\Phi}
\left(
j_{l_{1}}^{(m)}(a_{1})
\ldots 
j_{l_{2k}}^{(m)}(a_{2k})
\right)
=
\sum_{(m_{1}, \ldots , m_{2k})\in \Upsilon_{2k}^{m}}
\widetilde{\Phi}
\left(
j_{l_{1},m_{1}}(a_{1})\ldots j_{l_{2k},m_{2k}}(a_{2k})
\right) .
$$
Keeping in mind that
$j_{l_{i},m_{i}}(a_{i})=
\gamma_{m_{i}}^{(l_{i})}(a)-\widehat{\gamma}_{m_{i}}^{(l_{i})}(a)$,
$1\leq i \leq 2k$, we can see that the only way to avoid a
separation of $a_{1}$ from
$a_{2k}$ (which would produce two singletons and thus give zero
contribution) is to take into account in the above sum only those tuples 
$(m_{1}, \ldots , m_{2k})\in \Upsilon_{2k}^{m}$, for which
$m_{2}, \ldots , m_{2k-1}\neq 1$ (i.e. in particular, 
$m_{2}=m_{2k-1}=2$), and moreover, assume
that the products start 
with $\gamma_{2}^{(l_2)}(a_{2})$ and end with 
$\gamma_{2}^{(l_{2k-1})}(a_{2k-1})$. Then, at site $(l_{1},1)$
we get $a_{1}a_{2k}$ and at $(l_{p}, 1)$, $p\in \{2,\ldots ,
k\}$, we get either the projection $t$ or the unit ${\bf 1}$ and
$\widetilde{\phi}$ sends them to $1$. Therefore, we obtain
$$
\widetilde{\Phi}
\left(
j^{(m)}_{l_{1}}(a_{1})
\ldots 
j^{(m)}_{l_{2k}}(a_{2k})
\right)
=
\phi(a_{1}a_{2k})\widetilde{\Phi}
(j^{(m-1)}(a_{2})
\ldots  j^{(m-1)} (a_{2k-1}))
$$
$$
=\phi(a_{1}a_{2k})\prod_{i=1}^{k-1}\phi(a_{P_{i}'})=
\prod_{i=1}^{k}\phi(a_{P_{i}})
$$
by the inductive assumption with respect to $m$,
where 
$$
P=\{P_{1}, \ldots , P_{k}\},\;\;\;
P'=\{P_{2}, \ldots , P_{k}\}
$$ 
and $P_{1}=\{1,2k\}$.
\hfill $\Box$\\[0pt]
\indent{\par}
{\sc Lemma 3.4.}
{\it Assume that the partition $P$ associated 
with the tuple $(l_{1}, \ldots , l_{n})$, where $n=2k$,
is a crossing pair-partition.
If $\phi(a_{i})=0$ for $i=1, \ldots , n$, then}
$$
\widetilde{\Phi}
\left(
j_{l_{1}}^{(m)}(a_{1})
\ldots 
j_{l_{2k}}^{(m)}(a_{n})
\right)
=0
$$
{\it Proof.}
We will show that the correlation which corresponds to a
crossing pair-partition $P$ of $\{1, \ldots , 2k\}$ 
produces a singleton and thus vanishes by the mean zero
assumption. 

There exist $1\leq p<q<r<s \leq 2k$ such that 
$l_{p}=l_{q}=l$, $l_{r}=l_{s}=l'$.
It is enough to consider those terms from the ``pyramid'' in
which $m_{p}=m_{q}=u$ and
$m_{r}=m_{s}=w$ since otherwise
we obtain at least one singleton which makes the contribution
vanish. Suppose now that $u\leq w$. Then
$j_{l,u}(a_{p})$ and
$j_{l,u}(a_{q})$ put a projection $t$ at
site $(l',w)$ since they put a $t$ at all sites $(b,c)$, where
$b\neq l$ and $ c\geq u$. Thus, at site $(l',w)$ we obtain
$$
\ldots t\ldots a_{r}\ldots t \ldots a_{s}\ldots 
$$
and thus $t$ separates $a_{r}$ and $a_{s}$. If $u>w$,
then a similar thing happens to $a_{p}$ and $a_{q}$ 
at site $(l,u)$. This makes the
contribution of all terms vanish.\hfill $\Box$

Assume now that ${\cal A}_{l}={\cal A}$, $l\in {\nz}$.
We will derive the central limit theorem for the sums of
$m$-free ``independent'' variables (in other words, the central
limit theorem for $m$-free convolutions) 
$$
S_{N}^{(m)}(a)=\frac{1}{\sqrt{N}}\sum_{k=1}^{N}j_{k}^{(m)}(a),
$$
where $a\in {\cal A}^{0}$.\\[0pt]
\indent{\par}
{\sc Theorem 3.5.}
{\it Let $m\in {\nz}$, $a_{1}, \ldots , a_{n} \in {\cal A}$, 
and let $\phi$ be a state on ${\cal A}$ for which $\phi(a_{i})=0$, 
$i=1, \ldots , n$. Then}
$$
\lim_{N\rightarrow \infty}\widetilde{\Phi} 
\left(
S_{N}^{(m)}(a_{1})
\ldots
S_{N}^{(m)}(a_{n})
\right)=
\sum_{\{P_{1}, \ldots , P_{k}\}
\in NC_{n}^{{\rm pair}}(m)}\phi(a_{P_{1}})\ldots \phi(a_{P_{k}})
$$
{\it if $n=2k$. If $n$ is odd, then the above limit vanishes.}\\[5pt]
{\it Proof.}
Using Proposition 2.2 and typical central limit arguments (see, 
for instance, the limit theorem for correlations which are
invariant under order-preserving injections in [Len98] or
[SvW94]) we know that only pair partitions may give a
nonvanishing contribution as $N\rightarrow \infty$.
Now use Lemmas 3.2-3.4 to see that out of these only
the non-crossing pair partitions of depth $\leq m$ 
really do give a nonvanishing contribution.
The second part of the theorem is again standard and
follows from the assumption on the zero mean.
\hfill $\Box$\\[0pt]
\indent{\par}
{\sc Corollary 3.6.}
{\it In particular, if ${\cal A}={\cz}[x]$, $x^{*}=x$, and
$\phi(x^{2})=1$, then}
$$
M_{n}^{(m)}\equiv \lim_{N\rightarrow \infty}\widetilde{\Phi}
\left(
(S_{N}^{(m)}(x))^{n}
\right)=|NC_{n}^{{\rm pair}}(m)|
$$
{\it for $n$ even. The odd limit moments vanish.}\\[5pt]
{\it Proof.} It follows immediately from Theorem 3.5.\\[0pt]
{\it Remark.}
Knowing that $m$-freeness approximates freeness, we
automatically obtain the central limit theorem for free random
variables (as well as conditionally free random variables or
their possible generalizations as discussed in [FLS98]). 
For that purpose and for given $n=2k$ it is enough to take the
$k$-free product state.

In Section 5 we will solve the moment problem for the limit moments
given by Corollary 3.6 for each $m$.\\[5pt]
\begin{center}
{\sc 4. An Invariance Principle and $m$-Free Brownian Motions}
\end{center}
In this section we state an invariance principle for the hierarchy
of freeness. We also define a corresponding {\it hierarchy of Brownian
motions} and show that under some additional assumptions on the state 
$\phi$, the limit distribution obtained from the invariance principle
are the distributions of the hierarchy of Brownian motions.

Let us begin with the invariance principle. Let $a\in {\cal
A}^{0}$ and instead of the sums $S_{N}^{(m)}(a)$, consider now
sample sums 
$$
S_{N,f}^{(m)}(a)=
\frac{1}{\sqrt{N}}
\sum_{k=1}^{\infty}
j^{(m)}_{k}(a) \;\int_{k-1}^{k}f
\left(\frac{t}{N}\right)dt,
$$
indexed not only by $N$ and $m$, but also by 
$f\in L^{2}_{c}(\rz _{+})$, where
$L^{2}_{c}(\rz_+)$ stands for the square integrable 
real-valued functions with compact support on $\rz$.\\[0pt]
\indent{\par}
{\sc Theorem 4.1.}
{\it Let $f_{1}, \ldots , f_{n}\in L^{2}_{c}(\rz_{+})$,
$a_{1}, \ldots ,a_{n}\in {\cal A}^{0}$, $m,N \in \nz$. Then}
$$
\lim_{N\rightarrow \infty}
\widetilde{\Phi}
\left(
S_{N,f_{1}}^{(m)}(a_{1})\ldots S_{N,f_{n}}^{(m)}(a_{n})
\right)
$$
$$
=
\sum_{\{P_{1}, \ldots , P_{k}\}\in NC_{n}^{{\rm pair}}(m)}
\phi(a_{P_{1}})\ldots \phi(a_{P_{k}})
\prod_{r=1}^{k}
\int_{0}^{\infty}
f_{\alpha(r)}(t)f_{\beta(r)}(t)dt
$$
{\it if $n=2k$, where $P_{i}=\{\alpha(i), \beta(i)\}$, $i=1, \ldots ,
k$. If $n$ is odd, then the above limit vanishes.}\\[5pt]
{\it Proof.}
This is a special case of the invariance principle for 
correlations invariant under order preserving injections proved
in [SvW94]. 
\hfill $\Box$

Under certain additional assumptions one can realize 
the limit distribution in terms of creation and annihilation 
operators on a suitable Fock space. Note that the only
difference between our invariance principle and the invariance
principle for free independence is that in the case of
$m$-freeness only non-crossing partitions of depth $\leq m$
survive in the limit. 

To take that into account it is enough
to define the $m$-{\it free Fock space}
$$
{\cal F}^{(m)}\equiv {\cal F}^{(m)}(L^2(\rz_+))=\cz 
\oplus \bigoplus_{k=1}^m
L^2(\rz_+)^{\otimes k}
$$
with the vacuum vector $\Omega_{m}=1\oplus 0 \oplus \ldots \oplus 0$
and the canonical scalar product $\langle .,.\rangle _{{\cal F}^{(m)}}$.

Next, we define the $m$-{\it free creation operators}
$$
a^{(m)*}(f):\; {\cal F}^{(m)}\rightarrow {\cal
F}^{(m)}
$$
$$
a^{(m)*}(f)\; f_{1}\otimes \ldots \otimes f_{n} =
\left\{
\begin{array}{cll}
f\otimes f_{1}\otimes \ldots \otimes f_{n}& {\rm if} & 1\geq n < m\\
0& {\rm if} & n=m
\end{array}
\right.
$$
with $a^{(m)*}(f)\Omega_{m}=f$ and the $m$-{\it free
annihilation operators}
$$
a^{(m)}(f):\; {\cal F}^{(m)}\rightarrow {\cal
F}^{(m)}
$$
$$
a^{(m)}(f)\; f_{1} \otimes \ldots \otimes f_{n} =
\langle f , f_{1} \rangle \; f_{2} \otimes \ldots \otimes f_{n}
$$
if $ 1 \leq n \leq m$ and $a^{(m)}(f)\Omega =0$.
Note that $a^{(m)*}(f), a^{(m)}(f)\in B({\cal F}^{(m)})$.

We are ready to find a realization of the invariance principle 
limit in terms of the $m$-free creation and annihilation
operators under standard assumptions.
For simplicity we assume that ${\cal A}$ is the *-algebra
generated by one element $a$, which we denote 
${\cal A} =\cz \langle a,a^{*} \rangle$.\\[0pt]
\indent{\par}
{\sc Theorem 4.2.}
{\it Let $\phi$ be a state on $\cz\langle a,a^*\rangle$ such that 
$\phi(a)=\phi(a^*)=\phi(aa)=\phi(a^*a)=\phi(a^*a^*)=0$,
$\phi(aa^*)=1$. Then,}
$$
\lim_{N\to\infty} \widetilde{\Phi} 
\left( 
S_{N,f_{1}}^{(m)}(a^{\varepsilon_1})\ldots
S_{N,f_{n}}^{(m)}(a^{\varepsilon_n})
\right)
= 
\langle \Omega_{m}, a^{(m)\varepsilon_1}(f_{1}) \ldots
a^{(m)\varepsilon_n}(f_{n})\Omega_{m} \rangle_{{\cal F}^{(m)}}
$$
{\it for all $n\in\nz$, $a^{\varepsilon_1},
\ldots,a^{\varepsilon_n}\in\{a,a^*\}$, 
$f_{1}, \ldots , f_{n}\in L^{2}_{c}(\rz _{+})$.}\\[5pt]
{\it Proof:} It is enough to notice that the $m$-truncated
creation and annihilation operators are defined in such a way that
there can be no contribution from pair-partitions of depth greater than $m$
since the latter would require a tensor product of order 
greater than $m$. 
\hfill $\Box$

For each $m\in \nz$ 
denote by ${\cal C}^{(m)}$ the $C^{*}$-algebra generated by 
$a^{(m)*}(f), a^{(m)}(f)$, $f\in L^{2}(\rz _{+})$ and let 
$\varphi_{m}$ be the vacuum expectation state in the $m$-free
Fock space. Then the pair $({\cal C}^{(m)}, \varphi_{m})$ 
can be viewed as the $m$-{\it free Brownian motion} and the
collection $({\cal C}^{(m)}, \varphi_{m})_{m\in \nz}$ as the 
{\it hierarchy of $m$-free Brownian motions}. \\[5pt]
\begin{center}
{\sc 5. The Hierarchy of Limit Measures}
\end{center}
In this section we solve the moment problem for the limit laws
obtained in the central limit theorem in the case when 
${\cal A}=\cz [x]$, where $x=x^{*}$. 
We obtain a sequence $(\mu_{m})_{m \in \nz}$
of discrete measures that approximate the Wigner measure.

For that purpose, let us introduce the 
{\it hierarchy of Cauchy transforms} $(G_{m}(z))_{m \in \nz}$
for the sequence of limit laws given by Corollary 3.6:
$$
G_m(z) =\sum_{n=0}^{\infty}M^{(m)}_{n}z^{-n-1}
$$
where $M^{(m)}_{n}=|NC_{n}^{{\rm pair}}(m)|$, $m, n\in \nz$, and,
in addition $M_{0}^{(m)}=1$, $m\in \nz $.
We also adopt the convention that $M_{n}^{(0)}=\delta_{n,0}$
which gives $G_{0}=1/z$. For the use of Cauchy transforms in the
case of freeness (conditional freeness), see [Voi86] and [Maa92]
([BLS96]).

The moments $M_{n}^{(m)}$ grow less
rapidly as $N\rightarrow \infty$ than the moments $M_{n}$ of the
Wigner measure, therefore it is clear that for each $m$ there
exists a unique measure $\mu^{(m)}$ of which $G_{m}$ is the
Cauchy transform. In particular, $\mu^{(0)}=\delta_{0}$.
We will find the explicit form of $\mu^{(m)}$ for each $m\in \nz$.
\\[0pt]
\indent{\par}
{\sc Lemma 5.1.}
{\it The hierarchy of Cauchy transforms satisfies the recurrence
relation} 
$$
G_{m}(z) = \frac{1}{z - G_{m-1}(z)},
$$
{\it where $m \in \nz$, with $G_{0}(z)=1/z$, if $\;{\rm Im}z\neq 0$.}\\[5pt]
{\it Proof:}
Let us assume that we know the number of non-crossing pair
partitions of depth less than or equal to 
$m$ of the set $\{1,\ldots, 2k\}$ for any $k\le n$. To get a
non-crossing pair partition of depth less than or equal to $m$
of the set $\{1,\ldots,2n+2\}$, we have to 
choose a number $k\in\{2,\ldots,2n+2\}$ that will form a pair with
$1$, then choose a non-crossing pair partition of depth less than
or equal to $m-1$ for the numbers between $1$ and $k$, i.e. of the set
$\{2,\ldots,k-1\}$, and a non-crossing pair partition of depth less than
or equal to $m$ for the numbers from $k+1$ to $2n+2$, i.e. of the
set $\{k+1,\ldots,2n+2\}$. 

Therefore, there are exactly
$|NC_{k-2}^{{\rm pair}}(m-1)| \, |NC_{2n-k+2}^{{\rm pair}}(m)|$ 
such pair partitions in which $1$ is paired with $k$. 
For the total number of non-crossing pair partitions
of depth less than or equal to $m$ of the set
$\{1,\ldots,2n+2\}$ we get
$$
|NC_{2n+2}^{{\rm pair}}(m)| = \sum_{k=2}^{2n+2} |NC_{k-2}^{{\rm
    pair}}(m-1)| \, |NC_{2n-k+2}^{{\rm pair}}(m)|.
$$
The terms with odd $k$ give zero since there can be no pair
partition of a set with an odd number of elements. Hence,
$$
M^{(m)}_{2n+2}=\sum_{k=2}^{2n+2}M^{(m-1)}_{k-2}M^{(m)}_{2n-k+2}=
\sum_{l=1}^{n+1}M^{(m-1)}_{2l-2}M^{(m)}_{2n-2l+2}.
$$
The recurrence relation for the moments leads easily to
the desired recurrence relation for the Cauchy transforms if
${\rm Im}z \neq 0 $ for since
\begin{eqnarray*}
G_{m}(z) &=& \sum_{n=0}^\infty M^{(m)}_{2n} z^{-2n-1}  =
\frac{1}{z} 
+ \sum_{n=0}^\infty M^{(m)}_{2n+2} z^{-2n-3}\\
&=& 
\frac{1}{z} + \frac{1}{z} \sum_{n=0}^\infty \sum_{l=1}^{n+1}
M^{(m-1)}_{2l-2} z^{-2l+1} \, 
M^{(m)}_{2n-2l+2} z^{-2n+2l-3} \\
&=&\frac{1}{z}+ \frac{G_m(z) G_{m-1}(z)}{z}
\end{eqnarray*}
and therefore
\[
G_m(z) = \frac{1/z}{1-G_{m-1}(z)/z} = \frac{1}{z-G_{m-1}(z)},
\]
which finishes the proof.
\hfill $\Box$\\
\indent{\par}
{\sc Remark 1.} Note that the series given by $G_{m}(z)$ converges
absolutely for $|z|> 2$ and all $m\in\nz$ since 
$$
|NC_{2k}^{\rm pair}(m)|\leq |NC_{2k}^{\rm pair}|,
$$
where 
$$
|NC_{2k}^{\rm pair}|=\frac{1}{k+1}{2k \choose k}
$$
denotes the number of {\it all} non-crossing partitions
of the set $\{1, \ldots , 2k\}$. Clearly, $|NC_{n}^{\rm pair}|=
|NC_{n}^{\rm pair}(m)|=0$ if $n$ is odd.\\[0pt]
\indent{\par}
{\sc Remark 2.}
The Cauchy transforms $G_{m}(z)$ are rational functions of the
complex variable $z$. In particular, 
$$
G_{0}(z)=\frac{1}{z}, \;\; G_{1}(z)=\frac{1}{z-\frac{1}{z}},\;\;
G_{2}(z)=\frac{1}{z-\frac{1}{z-\frac{1}{z}}}, \ldots  .
$$
\indent{\par}
We will show below that $G_m$ has $m+1$ simple poles in the
interval $(-2,2)$ (and none anywhere else). 
For that purpose we use the Chebyshev polynomials of the
second kind 
$$
U_m(x)=\frac{\sin[(m+1)\arccos(x)]}{\sin(\arccos(x))},
$$
for $x\in (-1,1)$, $m\in \nz \cup \{0\}$.
They satisfy the recurrence relation
$$
U_{m+1}(x) = 2xU_m(x) - U_{m-1}(x)
$$
with $U_0(x)=1$. Denote by $U_m(z)$ the analytic extension of
$U_m(x)$. Note that $U_m(z)$ has exactly $m$ simple zeros 
$$
u_{m,k}=\cos\left(\frac{k \pi}{m+1}\right),\;\;\; k=1,\ldots, m 
$$
and that the zeros of $U_{m}(z)$ differ from those of
$U_{m+1}(z)$. This enables us to define the meromorphic
function
$$
W_{m}(z)=\frac{U_{m}(z/2)}{U_{m+1}(z/2)},\;\;\; 
m \in \nz \cup \{0\},
$$
with $m+1$ simples poles on the real line given by 
$$
z_{m,k}=2 \cos\left(\frac{k \pi}{m+2}\right),\;\;\;
k=1,\ldots,m+1.
$$
We show below that $W_{m}(z)$ coincides with $G_{m}(z)$.\\[0pt]
\indent{\par}
{\sc Lemma 5.2.}
{\it Let $m\in \nz \cup \{0\}$. 
The Cauchy transform $G_{m}(z)$ agrees with $W_{m}(z)$ for
$z\notin \{z_{m,k}|\; 1\leq k \leq m+1\}$.}\\[5pt]
{\it Proof:} 
Clearly, $W_{0}(z)=G_{0}(z)=1/z$ 
since $U_{0}(z)=1$ and $U_{1}(z)=2z$.
Let us show that the functions $W_{m}(z)$
satisfy the recurrence relation given by Lemma 5.1. 
If $m\geq 1$, then the recurrence relation for the Chebyshev
polynomials of the second kind gives
\begin{eqnarray*}
W_{m+1}(z)&=& 
\frac{U_{m+1}(z/2)}{U_{m+2}(z/2)}
=
\frac{U_{m+1}(z/2)}{z U_{m+1}(z/2) - U_m(z/2)} \\
&=& \frac{1}{z-U_m(z/2)/U_{m+1}(z/2)}
= \frac{1}{z-W_{m}(z)}
\end{eqnarray*}
for all $z\notin \{z_{m,k}|\; 1\leq k \leq m+1\}$.
Therefore, $G_{m}(z)$ must agree with 
$W_{m}(z)$ also for $m\geq 1$ on the intersection of their
domains, therefore, by uniquness of analytic continuation, they
must have the same domain, which finishes the proof. 
\hfill $\Box$\\[0pt]
\indent{\par}
{\sc Theorem 5.3}
{\it The measures $\mu^{(m)}$ take the form} 
$$
\mu^{(m)} = \sum_{k=1}^{m+1} b_{m,k} \delta_{z_{m,k}}
$$
{\it where}
$$
b_{m,k}=\frac{2 \sin^2\left(k\pi/(m+2)\right)}{m+2}
$$ 
{\it for $m\in \nz \cup \{0\}$ and $k=1,\ldots,m+1$}.\\[5pt]
{\it Proof:} 
We have to invert the Cauchy tranforms. 
By Lemma 5.2, $G_{m}(z)$ is a 
rational function with the degree of the denominator exceeding that
of the numerator and with simple poles at $z_{m,k}$, $1\leq k
\leq m+1$. 
Thus its decomposition into partial fractions takes the form
$$
G_m(z) = \sum_{k=1}^{m+1} \frac{b_{m,k}}{z-z_{m,k}}.
$$
This shows that $G_{m}(z)$ is the Cauchy transform of a discrete
measure with point masses at $z_{m,k}$, $1\leq k \leq m+1$. 
The calculation of the residues gives the masses
\begin{eqnarray*}
b_{m,k}
&=&
\lim_{z\rightarrow z_{m,k}}
\frac
{\sin [(m+1)\arccos (z/2)]}
{d/dz\sin [(m+2)\arccos (z/2)]}\\
&=&
\frac{2 \sin^2\left(k\pi/({m+2})\right)}{m+2}.
\end{eqnarray*}
which finishes the proof.
\hfill $\Box$\\[0pt]
\indent{\par}
{\sc Example.} The measures $\mu^{(0)}, \mu^{(1)}, \mu^{(2)}$
are given by
$$
\mu^{(0)}=\delta_{0},\;\; \mu^{(1)}=\frac{1}{2}\delta_{-1} +
\frac{1}{2}\delta_{1},\;\;
\mu^{(2)}=\frac{1}{4}\delta_{-\sqrt{2}} +\frac{1}{2}\delta_{0}+
\frac{1}{4}\delta_{\sqrt{2}}.
$$
\indent{\par}
Since the moment problems are determined for all $m\in \nz$,
i.e. the measures $\mu^{(m)}$ are uniquely determined,
$\mu^{(m)}$ converges weakly to the Wigner measure $\mu_{W}$.\\[5pt]
\begin{center}
{\sc 6. Poisson's Limit Theorem}
\end{center}
In this section we study Poisson's limit theorem for the
hierarchy of freeness and solve the moment problems for the
limit laws. By $|NC_{n}(b,m)|$ we denote the number of non-crossing
partitions of $\{1,\ldots,n\}$ with $b$ blocks and depth less than or
equal to $m$.\\[0pt]
\indent{\par}
{\bf Theorem 6.1.}
{\it Let ${\cal A}_{l}={\cal A}=\cz [x]$, $l\in L$, $x^{*}=x$, 
and assume that $N\phi^{N}(a^{k})\rightarrow \lambda$, $k\in
\nz$, $\lambda >0$. Let $S_{m,N}=\sum_{k=1}^{N}j_{k}^{(m)}(x)$
and denote by $\widetilde{\Phi}^{(m,N)}$ the $m$-free product
state corresponding to $\phi^{N}$. Then}
$$
\lim_{N\rightarrow \infty}\widetilde{\Phi}^{(m,N)}(S_{m,N}^{n})=
\sum_{q=1}^{n}\lambda^{q}|NC_{n}(q,m)|\equiv M_{n}^{(m)}(\lambda)
$$
{\it Proof.}
We have
$$
\widetilde{\Phi}^{(m,N)}(S_{m,N}^{n})=
\sum_{1\leq k_{1}, \ldots , k_{n}\leq N}
\widetilde{\Phi}^{(m,N)}(j_{k_{1}}(x)\ldots j_{k_{n}}(a))
=\sum_{P\in {\cal P}_{n}}(N)_{b(P)}m(P)
$$
where $P_{n}$ denotes partitions of $\{1, \ldots , n\}$,
$m(P)=\widetilde{\Phi}^{(m,N)}(j_{k_{1}}(x)\ldots j_{k_{n}}(x))$
for any tuple $(k_{1}, \ldots , k_{n})$ 
associated with the partition $P$, $b(P)$ denotes the number of
blocks of $P$ and $(N)_{r}=N(N-1)\ldots (N-r+1)$.

Now we apply the usual Poisson's limit arguments. 
The only partitions $P$ which survive in the limit $N\rightarrow
\infty$ are those for which the expression for $m(P)$ contains
a term of type $\lambda^{b(P)}$ (i.e. the number of blocks of $P$ is
equal to the number of moments in the given term).
If $P$ is a crossing partition then $m(P)$ ``factorizes'' into more than
$b$ moments and thus gives no contribution to the limit. If $P$
is non-crossing, then we have two cases: (i) $d(P)>m$ and (ii)
$d(P)\leq m$. In case (i) the contribution is zero even before
taking the limit by the GNS construction. In case (ii) the
contribution is $\lambda^{b(P)}$, which ends the proof.
\hfill $\Box$

In order to solve the associated moment problem, we 
want to find the generating functions for
$|NC_{n}(b,m)|$. Thus, let 
$$
H^{(m)}(\lambda,z)=\sum_{n,b=0}^\infty
|NC_n(b,m)| \lambda^b z^{-n-1}
$$
for $m\geq 1 $ and $H^{(0)}(\lambda,z)=1/z$, where we adopt the conventions
that $|NC_n(b,0)|=\delta_{n0}\delta_{b0}$ and $|NC_n(0,m)| = \delta_{n0}$.
Clearly $|NC_n(b,m)|=0$ for $b>n>0$, so the summation over $b$
is finite for fixed $n$.

Note that $H^{(m)}(z)$, $m\geq 0$, converge absolutely for $|z|$
sufficiently large, say $|z|>R(\lambda)=(\sqrt{\lambda}+1)^2$.
Moreover, they go to zero as $|z|$ goes to infinity
(since there is no constant term in the series). Thus
$|H^{(m)}(\lambda,z)| < 1$ for $|z|> R'(\lambda)$ for some sufficiently
large $R'(\lambda)$ (it depends on $\lambda$ but not on $m$ by
comparison with the free Poisson law, i.e.
$|NC_n(b,m)| \le |NC_n(b)|$ and therefore $|H^{(m)}(\lambda,z)|\le
H(|\lambda|,|z|)$, where $|NC_{n}(b)|$ denotes the number
of non-crossing partitions of $\{1, \ldots , n\}$ of $b$ blocks
and $H(\lambda, z)$ is the generating function for the free
Poisson law).\\[0pt] 
\indent{\par}
{\sc Lemma 6.2.}
{\it The hierarchy of generating functions $(H^{(m)})_{m\geq 0}$
satisfies the recurrence relation}
$$
H^{(m)}(\lambda,z) = \frac{1-H^{(m-1)}(\lambda,z)}{z-zH^{(m-1)}(\lambda,z) -
  \lambda}
$$
{\it for $m=1,2,\ldots$ and $|z|>R'(\lambda)$.}\\[0pt]

{\it Proof:} To get a non-crossing partition of $\{1,\ldots,n\}$
($n\ge 1$) we pick the elements that will be put in the same block as
the first element, denote this block by $\{1,1+k_1,1+k_1+k_2,\ldots
1+k_1+\cdots+k_{r-1}\}$, and then choose non-crossing partitions for
the remaining intervals $\{2,\ldots,k_1\}$, $\{k_1+2, \ldots,
k_1+k_2\}$, $\ldots$,
$\{k_1+\cdots+k_{r-2}+2,\ldots,k_1+\cdots+k_{r-1}\}$,
$\{k_1+\cdots+k_{r-1}+2,\ldots, n\}$. We will denote the number of
elements of the last interval by $k_r$. If we want the resulting
partition to have depth $\le m$, then the partitions chosen for  $\{2,\ldots,k_1\}$, $\ldots$,
$\{k_1+\cdots+k_{r-2}+2,\ldots,k_1+\cdots+k_{r-1}\}$ must have depth
$\le m-1$, and that chosen for $\{k_1+\cdots+k_{r-1}+2,\ldots, n\}$
must have depth $\le m$. Let $b_k$ be the number of blocks of the
partition of the $k^{{\rm th}}$ interval, then the number of blocks of
the whole partition is $b_1+\cdots+b_r+1$. Therefore the number of
non-crossing partitions of $\{1,\ldots,n\}$ with $b$ blocks and depth
$\le m$ can be calculated recursively by the formula
$$
|NC_n(b,m)| =
\sum_{r=1}^n 
\sum_{\stackrel{\scriptstyle k_1,\ldots,k_{r-1}\ge 1; k_r\ge 0}
{k_1+\cdots+k_r=n-1}}
\sum_{\stackrel{\scriptstyle b_1,\ldots,b_r\ge 0}{b_1+\cdots+b_r=b-1}} 
$$
$$
|NC_{k_1-1}(b_1,m-1)|\cdots
|NC_{k_{r-1}-1}(b_{r-1},m-1)|\, |NC_{k_r}(b_r,m)| 
$$
for $n\ge 1$, if we use the conventions
$|NC_n(b,0)|=\delta_{n0}\delta_{b0}$ and $|NC_n(0,m)| = \delta_{n0}$.
By these conventions we have $H^{(0)}(\lambda,z)=1/z$.

Let now $m\ge 1$ and $|z| \ge R'$. Then we have
\begin{eqnarray*}
H^{(m)}(\lambda,z) &=& \sum_{n,b=0}^\infty |NC_n(b,m)| \lambda^b z^{-n-1}
\\
&=& \frac{1}{z} +  \frac{\lambda}{z}\sum_{n,b=1}^\infty \sum_{r=1}^n
\sum_{\stackrel{\scriptstyle k_1,\ldots,k_{r-1}\ge 1; k_r\ge
    0}{k_1+\cdots+k_r=n-1}} \sum_{\stackrel{\scriptstyle 
b_1,\ldots,b_r\ge 0}{b_1+\cdots+b_r=b-1}} |NC_{k_1-1}(b_1,m-1)|
\lambda^{b_1} z^{-k_1} \times \\
&& \quad \cdots
|NC_{k_{r-1}-1}(b_{r-1},m-1)|\lambda^{b_{r-1}} z^{-k_{r-1}}
|NC_{k_r}(b_r,m)| \lambda^{b_r} z^{-k_r-1} \\
&=& \frac{1}{z} + 
\frac{\lambda}{z} \sum_{r=1}^\infty 
\left( 
\sum_{\beta,\nu=0}^\infty
|NC_\nu(\beta,m-1)| \lambda^\beta z^{-\nu-1} 
\right)^{r-1} 
\sum_{\mu,\alpha=0}^\infty
|NC_{\mu}(\alpha,m)| \lambda^{\alpha} z^{-\mu-1} \\
&=& \frac{1}{z} + 
\frac{\lambda H^{(m)}(\lambda,z)}{z\big(1-H^{(m-1)}(\lambda,z)\big)},
\end{eqnarray*}
where the summations can be interchanged since all sums converge
absolutely (remember that 
$|H^{(m-1)}(\lambda,z)| < 1$ for $|z| >R'(\lambda)$), 
and therefore
$$
H^{(m)}(\lambda,z) = \frac{1-H^{(m-1)}(\lambda,z)}{z - zH^{(m-1)}(\lambda,z)
-\lambda}.
$$
\hfill $\Box$

We will now give an explicit expression for the solution of this
recurrence relation. To this end we will again use the Chebyshev
polynomials of the second kind. \\[0pt]
\indent{\par}
{\sc Proposition 6.3.}
{\it Let $\lambda>0$, $m\in \nz\cup\{0\}$. The meromorphic functions}
$$
F_\lambda^{(m)}(z) =
\frac{(z-\lambda)U_{m}\left(\frac{z-\lambda-1}{2\sqrt{\lambda}}\right)
  - \sqrt{\lambda}
  U_{m+1}\left(\frac{z-\lambda-1}{2\sqrt{\lambda}}\right)}{
  z U_m\left(\frac{z-\lambda-1}{2\sqrt{\lambda}}\right)}
$$
{\it solve the recurrence relation $F_\lambda^{(m)}(z) =
\frac{1-F_\lambda^{(m-1)}(z)}{z-zF_\lambda^{(m-1)}(z) -\lambda}$, for $m\ge
1$, $F_\lambda^{(0)}(z) = 1/z$, and therefore we have
$H^{(m)}(\lambda,z) = F^{(m)}_\lambda(z)$ for $|z|>R'(\lambda)$.

Furthermore, $F^{(m)}_\lambda(z)$ has the partial fraction decomposition
$$
F^{(m)}_{\lambda}(z) = 
\sum_{k=0}^m \frac{a_{m,k}(\lambda)}{z-y_{m,k}(\lambda)},
$$
where
\begin{eqnarray*}
y_{m,0}(\lambda) &=& 0, \\
y_{m,k}(\lambda) &=& 2\sqrt{\lambda} \cos\left(\frac{k\pi}{m+1}\right)+\lambda+1,
\qquad k=1,\ldots, m, \\
a_{m,0}(\lambda) &=& \sqrt{\lambda}
\frac{U_{m+1}\left( \frac{\lambda
      +1}{2\sqrt{\lambda}}\right)}{U_{m}\left( \frac{\lambda
      +1}{2\sqrt{\lambda}}\right)} - \lambda, \\
a_{m,k}(\lambda) &=&   \frac{2\lambda \sin^2
  \left(\frac{k\pi}{m+1}\right)}{(m+1)\left[2\sqrt{\lambda}\cos\left(\frac{k\pi}{m+1}\right) +\lambda +1\right]}, \qquad k=1,\ldots, m, \\
\end{eqnarray*}
for $m\in\nz$.}\\[5pt]
{\it Proof:}
Fix $\lambda$  and let 
$F_{\lambda}^{(m)}(z)=P_{\lambda}^{(m)}(z)/Q_{\lambda}^{(m)}(z)$,
where
\begin{eqnarray*}
  P_{\lambda}^{(m)}(z) &=& \lambda^\frac{m}{2} (z-\lambda) U_m\left(\frac{z-\lambda -
      1}{2\sqrt{\lambda}}\right) - \lambda^\frac{m+1}{2} U_{m+1}\left(\frac{z-\lambda -
      1}{2\sqrt{\lambda}}\right), \\
Q_{\lambda}^{(m)}(z) &=& \lambda^\frac{m}{2} z U_m\left(\frac{z-\lambda -
      1}{2\sqrt{\lambda}}\right).
\end{eqnarray*}
From the recurrence relation for the Chebyshev polynomials of the second
kind it follows that $P_{\lambda}^{(m)}(z), Q_{\lambda}^{(m)}(z)$ 
satify the coupled recurrence relations
\begin{eqnarray*}
P_{\lambda}^{(m)}(z) &=& Q_{\lambda}^{(m-1)}(z) - P_{\lambda}^{(m-1)}(z), \\
Q_{\lambda}^{(m)}(z) &=& (z-\lambda)Q_{\lambda}^{(m-1)}(z) - 
z P_{\lambda}^{(m-1)}(z),
\end{eqnarray*}
for $m\ge 1$, and $P_{\lambda}^{(0)}(z)=1$, $Q_{\lambda}^{(0)}(z) = z$.

For $m=0$ we have $F_{\lambda}^0(z) =
P_{\lambda}^{(0)}(z)/Q_{\lambda}^{(0)}(z) = 1/z$, and for $m\ge 1$
\begin{eqnarray*}
F^{(m)}_{\lambda}(z) &=& \frac{P_{\lambda}^{(m)}(z)}
{Q_{\lambda}^{(m)}(z)} = \frac{ Q_{\lambda}^{(m-1)}(z) -
  P_{\lambda}^{(m-1)}(z)}{(z-\lambda)Q_{\lambda}^{(m-1)}(z) - 
z P_{\lambda}^{(m-1)}(z)} \\
&=& \frac{1-P_{\lambda}^{(m-1)}(z)/Q_{\lambda}^{(m-1)}(z)}
{z-\lambda - z P_{\lambda}^{(m-1)}(z)/Q_{\lambda}^{(m-1)}(z)} =
\frac{1-F^{(m-1)}_{\lambda}(z)}{z-zF^{(m-1)}_{\lambda}(z)-\lambda}.
\end{eqnarray*}

It is easy to deduce from the recurrence relation that 
$P_{\lambda}^{(m)}(z)$ has
degree $\le m$. From the definition of 
$Q_{\lambda}^{(m)}(z)$ we immediately see
that is has $m+1$ distinct simple roots, $y_{m,0}(\lambda)=0$, and
$y_{m,k}(\lambda)=2\sqrt{\lambda}\cos\left(\frac{k\pi}{m+1}\right) +\lambda +
1$, $k=1,\ldots,m$. Therefore $F^{(m)}_{\lambda}(z)$ 
has the form stated in the
proposition. The calculation of the residues gives
\begin{eqnarray*}
a_{m,0}(\lambda) &=& \lim_{z\to 0} z F_\lambda^{(m)}(z) = \sqrt{\lambda}
\frac{U_{m+1}\left( \frac{\lambda
      +1}{2\sqrt{\lambda}}\right)}{U_{m}\left( \frac{\lambda
      +1}{2\sqrt{\lambda}}\right)} - \lambda, \\
a_{m,k}(\lambda) &=& \lim_{z\to y_{m,k}} (z-z_{m,k}) F_\lambda^{(m)}(z) \\
&=& - \frac{2\lambda}{2\sqrt{\lambda} \cos\left(\frac{k\pi}{m+1}\right)
  + \lambda +1} \lim_{x\to x_{m,k}}
\frac{\sin[(m+2)\arccos(x)]}{\frac{{\rm d}}{{\rm d} x}
  \sin[(m+1)\arccos(x)]} \\
&=& \frac{2\lambda \sin^2
  \left(\frac{k\pi}{m+1}\right)}{(m+1)\left[2\sqrt{\lambda}\cos\left(\frac{k\pi}{m+1}\right) +\lambda +1\right]}.
\end{eqnarray*}
for $m\geq k \geq 1$.
\hfill $\Box$\\[0pt]
\indent{\par}
{\sc Theorem 6.4.}
{\it Let $m\in\nz$, $\lambda>0$. 
The moments $\big(M^{(m)}_n(\lambda)\big)_{n\in\nz}$ determine a
unique measure on the real line of the form }
$$
\mu^{(m)}_{\lambda} = \sum_{k=0}^m a_{m,k}(\lambda) \delta_{y_{m,k}(\lambda)}.
$$
{\it Proof:}
The moments $(M^{(m)}_n(\lambda))_{n\in\nz}$ grow less rapidly as
$n\to\infty$ than the moments of the free Poisson limit measure,
therefore it is clear that the moment problem has a unique solution
$\mu^{(m)}_\lambda$. Denote its Cauchy transform by
$G^{(m)}_\lambda(z)=\int_\rz \frac{1}{z-x} {\rm d}\mu^{(m)}_\lambda(x)$.

By Lemma 6.2 we know that
$H^{(m)}(\lambda,z)=\sum_{n=0}^\infty M^{(m)}_n(\lambda) z^{-n-1}$
converges absolutely for $|z|\ge R(\lambda)$, therefore it coincides with the
Cauchy transform of $\mu^{(m)}_\lambda$ for $|z|\ge R(\lambda)$.
By Proposition 6.3 we now have $G^{(m)}_\lambda(z)=F^{(m)}_\lambda(z)$ for
$|z| > R'(\lambda)$, and then also for all $z\in\cz\backslash\rz$,
since both functions are analytic on 
$\cz\backslash\rz$. 

It now follows immediately from the partial fraction decomposition
of Proposition 6.3 that $\mu^{(m)}_\lambda$ has the form stated
in the theorem.
\hfill $\Box$\\[0pt]
\indent{\par}
{\sc Example:} We get
\begin{eqnarray*}
&&\mu^{(0)}_\lambda = \delta_0,\qquad
\mu^{(1)}_\lambda = \frac{1}{1+\lambda} \delta_0 +
\frac{\lambda}{1+\lambda} \delta_{1+\lambda}, \\
&&\mu^{(2)}_\lambda = \frac{1}{1+\lambda+\lambda^2}\delta_0 +
\frac{\lambda}{2(1+\sqrt{\lambda}+\lambda)}
\delta_{1+\sqrt{\lambda}+\lambda} +
\frac{\lambda}{2(1-\sqrt{\lambda}+\lambda)} \delta_{1-\sqrt{\lambda}
  +\lambda}.
\end{eqnarray*}
\vspace{20pt}
\begin{center}
{\sc References}
\end{center}
[BLS96] {\sc M.~Bo$\dot{{\rm z}}$ejko, M.~Leinert, R.~Speicher}, 
{\it Convolution and limit theorems for conditionally free
random variables}, Pac. J. Math. 175, No.2 (1996),
357-388.\\[3pt]
[FLS98] {\sc U.~Franz, R.~Lenczewski, M.~Sch\"{u}rmann},
{\it The GNS construction for the hierarchy of freeness}, Preprint No.
9/98, Technical University of Wroclaw, 1998.\\[3pt]
[Len97] {\sc R.~Lenczewski}, 
{\it Unification of independence in quantum probability},
Preprint No.\ 15/97, Technical University of Wroclaw, 1997, to appear in 
Inf.~Dim.~Anal.~Quant. Probab. \& Rel.~Top.\\[3pt]
[Len98] {\sc R.~Lenczewski}, 
{\it A noncommutative limit theorem for homogeneous
correlations}, Studia Math. 129 (1998),
225-252.\\[3pt] 
[Maa92] {\sc H.~Maassen}, 
{\it Addition of freely independent random variables},
J. Funct. Anal. 106 (1992), 409-438.\\[3pt]
[Sch94] {\sc M.~Sch\"{u}rmann}, {\it Non-commutative probability on
algebraic structures}, in "Probability measures on groups and
related structures", Vol. XI (Oberwolfach, 1994), 332-356, World
Scientific, River Edge, NJ, 1995.\\[3pt]
[Sch95] {\sc M.~Sch\"{u}rmann}, {\it Direct sums of tensor products and
non-commutative independence}, J. Funct. Anal. 133 (1995), 1-9.\\[3pt]
[Spe90] {\sc R.~Speicher}, 
{\it A new example of ``independence'' and ``white noise''},
Probab.~Th. Rel.~Fields 84 (1990) , 141-159.\\[3pt]
[SvW94] {\sc R.~Speicher, W.~von Waldenfels}, 
{\it A general central limit theorem and
invariance principle}, in: ``Quantum Probability and Related Topics'',
Vol.~IX, World Scientific, 1994, 371-387.\\[3pt]
[Voi85] {\sc D.V.~Voiculescu}, {\it Symmetries of some reduced free product 
${\cal C}^{*}$-algebras}, in ``Operator Algebras and their 
Connections with Topology and Ergodic Theory'', Lecture
Notes in Math. 1132, 556-588, Springer, Berlin, 1985.\\[3pt]
[Voi86] {\sc D.V.~Voiculescu}, {\it Addition of certain non-commuting
random variables}, J. Funct. Anal. 66 (1986), 323-346.
\end{document}